\UseRawInputEncoding

\documentclass[12pt,reqno]{amsart}
\usepackage{amsthm}

\usepackage{amssymb}
\usepackage{graphics}
\usepackage{tikz}
\usetikzlibrary{shapes,backgrounds,calc}
\usepackage{latexsym}
\usepackage{multicol}
\usepackage{verbatim,enumerate}
\usepackage{accents}
\usepackage{cite}

\usepackage{amsmath, amscd,url}

\usepackage{setspace}
\usepackage{pstricks}

\usepackage{longtable}            
\usepackage{rotating}
\usepackage{lscape}

\usepackage{xcolor,colortbl}
\definecolor{Gray}{gray}{0.85}
\newcolumntype{a}{>{\columncolor{Gray}}c}

\advance\textwidth by 1.3in \advance\oddsidemargin by -.6in \advance\evensidemargin by -.6in
\theoremstyle{definition}
\newtheorem{cor}{Corollary}
\newtheorem{lem}{Lemma}
\newtheorem{prop}{Proposition}
\newtheorem*{propn}{Proposition}
\newtheorem{thm}{Theorem}
\newtheorem*{theom}{Theorem}

\theoremstyle{definition}
\newtheorem{defn}{Definition}

\theoremstyle{definition}
\newtheorem{example}{Example}

\newtheorem{rem}{Remark}

\newenvironment{pf}{\proof}{\endproof}
\newcounter{cnt}
 \makeatletter
\def\mydggeometry{\makeatletter\dg@YGRID=1\dg@XGRID=20\unitlength=0.003pt\makeatother}
\makeatother \theoremstyle{remark}

\numberwithin{equation}{section}
\let\bwdg\bigwedge
\def\bigwedge{{\textstyle\bwdg}}

\newcommand{\nc}{\newcommand}
\newcommand{\rnc}{\renewcommand}

\nc{\cal}{\mathcal} \nc{\goth}{\mathfrak} \rnc{\bold}{\mathbf}

\nc\bomega{{\mbox{\boldmath $\omega$}}} \nc\bpsi{{\mbox{\boldmath $\Psi$}}}
 \nc\balpha{{\mbox{\boldmath $\alpha$}}}
 \nc\bpi{{\mbox{\boldmath $\pi$}}}
 \nc\bvpi{{\mbox{\boldmath $\varpi$}}}
\nc\chara{\operatorname{ch}}

  \nc\bxi{{\mbox{\boldmath $\xi$}}}
\nc\bmu{{\mbox{\boldmath $\mu$}}} \nc\bcN{{\mbox{\boldmath $\cal{N}$}}} \nc\bcm{{\mbox{\boldmath $\cal{M}$}}} \nc\blambda{{\mbox{\boldmath
$\lambda$}}}\nc\bnu{{\mbox{\boldmath $\nu$}}}

\makeatletter
\def\section{\def\@secnumfont{\mdseries}\@startsection{section}{1}%
  \z@{.7\linespacing\@plus\linespacing}{.5\linespacing}%
  {\normalfont\scshape\centering}}
\def\subsection{\def\@secnumfont{\bfseries}\@startsection{subsection}{2}%
  {\parindent}{.5\linespacing\@plus.7\linespacing}{-.5em}%
  {\normalfont\bfseries}}
\makeatother

 \nc{\Hom}{\operatorname{Hom}}
  \nc{\mode}{\operatorname{mod}}
\nc{\End}{\operatorname{End}} \nc{\wh}[1]{\widehat{#1}} \nc{\Ext}{\operatorname{Ext}} \nc{\ch}{\text{ch}} \nc{\ev}{\operatorname{ev}}
\nc{\Ob}{\operatorname{Ob}} \nc{\soc}{\operatorname{soc}} \nc{\rad}{\operatorname{rad}} \nc{\head}{\operatorname{head}}

\def\ann{\operatorname{Ann}}

 \nc{\Cal}{\cal} \nc{\Xp}[1]{X^+(#1)} \nc{\Xm}[1]{X^-(#1)}
\nc{\on}{\operatorname} \nc{\Z}{{\bold Z}} \nc{\J}{{\cal J}}  \nc{\Q}{{\bold Q}}

\nc{\N}{{\bold N}}  \nc\boa{\bold a} \nc\bob{\bold b} \nc\boc{\bold c} \nc\bod{\bold d} \nc\boe{\bold e} \nc\bof{\bold f} \nc\bog{\bold g}
\nc\boh{\bold h} \nc\boi{\bold i} \nc\boj{\bold j} \nc\bok{\bold k} \nc\bol{\bold l} \nc\bom{\bold m} \nc\bon{\mathbb n} \nc\boo{\bold o}
\nc\bop{\bold p} \nc\boq{\bold q} \nc\bor{\bold r} \nc\bos{\bold s} \nc\boT{\bold t} \nc\boF{\bold F} \nc\bou{\bold u} \nc\bov{\bold v}
\nc\bow{\bold w} \nc\boz{\bold z}\nc\ba{\bold A} \nc\bb{\bold B} \nc\bc{\mathbb C} \nc\bd{\bold D} \nc\be{\bold E} \nc\bg{\bold
G} \nc\bh{\bold H} \nc\bi{\bold I} \nc\bj{\bold J} \nc\bk{\bold K} \nc\bl{\bold L} \nc\bm{\bold M} \nc\bn{\mathbb N} \nc\bo{\bold O} \nc\bp{\bold
P} \nc\bq{\bold Q} \nc\br{\bold R} \nc\bs{\bold S} \nc\bt{\bold T} \nc\bu{\bold U} \nc\bv{\bold V} \nc\bw{\bold W} \nc\bz{\mathbb Z} \nc\bx{\bold
x} \nc\KR{\bold{KR}} \nc\rk{\bold{rk}} \nc\het{\text{ht }}

\nc\toa{\tilde a} \nc\tob{\tilde b} \nc\toc{\tilde c} \nc\tod{\tilde d} \nc\toe{\tilde e} \nc\tof{\tilde f} \nc\tog{\tilde g} \nc\toh{\tilde h}
\nc\toi{\tilde i} \nc\toj{\tilde j} \nc\tok{\tilde k} \nc\tol{\tilde l} \nc\tom{\tilde m} \nc\ton{\tilde n} \nc\too{\tilde o} \nc\toq{\tilde q}
\nc\tor{\tilde r} \nc\tos{\tilde s} \nc\toT{\tilde t} \nc\tou{\tilde u} \nc\tov{\tilde v} \nc\tow{\tilde w} \nc\toz{\tilde z} \nc\woi{w_{\omega_i}}

\newcommand\wrapped[1]%
{%
	\begin{array}{@{}l@{}}#1\end{array}%
}

\begin{document}


\title{A note on topological indices and the twin classes of graphs}

\author{P. Gangaeswari}
\address{Department of Mathematics, Manonmaniam Sundaranar University, Abishekapatti, Tirunelveli 627 012, Tamil Nadu, India.}
\email{amirdhavarshini11@gmail.com.}

\author{K. Selvakumar}
\address{Department of Mathematics, Manonmaniam Sundaranar University, Abishekapatti, Tirunelveli 627 012, Tamil Nadu, India.}
\email{selva\_158@yahoo.co.in.}

\author{G. Arunkumar$^{\ast}$}
\address{Department of Mathematics, Indian Institute of technology Madras, Tamilnadu, India.}
\email{arun.maths123@gmail.com, garunkumar@iitm.ac.in.}

\thanks{$^{\ast}$-The corresponding author. The first author acknowledges the CSIR Junior Research Fellowship (09/0652(11961)/2021-EMR-I). The third/corresponding author acknoledges the NFIG grant of Indian Institute of Technology Madras (RF/22-23/0985/MA/NFIG/009003​) }


\maketitle

\begin{abstract}
	Topological indices are parameters associated with graphs that have many applications in different areas such as mathematical chemistry.  Among various topological indices, the Wiener index is classical \cite{w}.  In this paper, we prove a formula for the Wiener index and more general $m$-Steiner Wiener index of an arbitrary graph $G$ in terms of the cardinalities of its twin classes. In particular, we will show that calculating these parameters for the graph $G$ can be reduced to calculating the same for a much smaller graph (in general) called the reduced graph of $G$. 

	
	
As applications of our main result, the $m$-Steiner Wiener index is explicitly calculated for various important classes of graphs from the literature  including

\begin{enumerate}
	\item[(a)] Power graphs associated with finite groups,
	\item[(b)] Zero divisor graphs and the ideal-based zero divisor graphs associated to commutative rings with unity, and
	\item[(c)] Comaximal ideal graphs associated with commutative rings with unity.
\end{enumerate}
We have also found an upper bound on the $m$-Steiner Wiener index of an infinite class of graphs called the completely joined graphs. As a corollary of this result, we explicitly calculate the $m$-Steiner Wiener index of the complete multipartite graphs.

\end{abstract}

\tableofcontents

\section{Introduction}

In graph theory, the notion of twin vertices plays an important role. In the literature, twin vertices have been studied extensively and proved to have applications in several branches \cite{iioa07, da08, ioa15, dp20}. Twin vertices are a useful tool in proving results in graph theory \cite{cmica10}.  For example, twin vertices naturally occur in the study of distance-hereditary graphs \cite{hh86}. Also, the twin free graphs are proved to be serving as a class of graphs with prescribed number negative eigenvalues in spectral graph theory \cite{gb21}. A generalization of twins can be seen in \cite{mrt14}.

We use twin vertices to write an arbitrary graph $G$ as a generalized composition of suitable choices of smaller graphs \cite{tfz08}. Let $G$ be an arbitrary graph. We define an equivalence relation $\sim^G$ on $V(G)$ as follows. For $u,v \in V(G)$, define $u \sim^G v$ if, and only if, $N_G(u) \backslash \{v\} = N_G(v) \backslash \{u\}$ where $N_G(u)$ is the open neighbourhood of the vertex $u$ in the graph $G$. In other words, two vertices of $G$ are related by $\sim^G$ if they are twin vertices in $G$. Let $C_1, \dots, C_k$ be the equivalence classes of this relation (with respective representatives $c_1,\dots,c_k$) which we call the twin classes of the graph $G$. 
\begin{propn} \cite[Lemma 3.1]{tfz08}
	Let $G$ be an arbitrary graph possibly with loops. Let the twin classes of the graph $G$ be $C_1,\dots,C_k$. Define $G_i := G[C_i]$ and $H := G[\{c_1,\dots,c_k\}]$. Then $G = H[G_1,\dots,G_k]$ where each $G_i$ is either complete or totally disconnected (empty) graph.
\end{propn}
The graph $H$ is said to be the reduced graph of the graph of $G$. 

\newpage
\begin{rem}
	In general, graph $H$ is much smaller than graph $G$ as illustrated in the following example. 

	
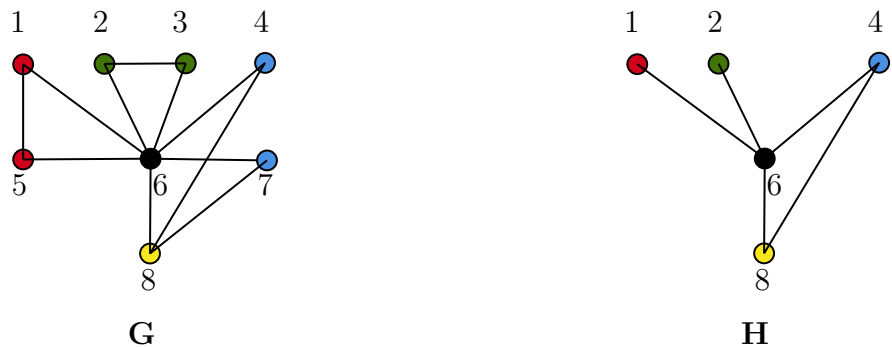
\begin{multicols}{2}
	\begin{center}
		\tikzset{every picture/.style={line width=0.75pt}} 
		
		\begin{tikzpicture}[x=0.75pt,y=0.75pt,yscale=-1,xscale=1]
			
			\draw  [fill={rgb, 255:red, 0; green, 0; blue, 0 }  ,fill opacity=1 ] (209.68,193.46) .. controls (206.92,193.48) and (204.66,191.25) .. (204.65,188.49) .. controls (204.63,185.73) and (206.86,183.48) .. (209.62,183.46) .. controls (212.38,183.44) and (214.63,185.67) .. (214.65,188.43) .. controls (214.66,191.19) and (212.44,193.44) .. (209.68,193.46) -- cycle ;
			\draw  [fill={rgb, 255:red, 74; green, 144; blue, 226 }  ,fill opacity=1 ] (267.39,145.12) .. controls (264.63,145.13) and (262.38,142.91) .. (262.36,140.15) .. controls (262.35,137.39) and (264.57,135.13) .. (267.33,135.12) .. controls (270.09,135.1) and (272.35,137.33) .. (272.36,140.09) .. controls (272.38,142.85) and (270.15,145.1) .. (267.39,145.12) -- cycle ;
			\draw  [fill={rgb, 255:red, 65; green, 117; blue, 5 }  ,fill opacity=1 ] (227.39,145.35) .. controls (224.63,145.37) and (222.38,143.15) .. (222.36,140.38) .. controls (222.35,137.62) and (224.57,135.37) .. (227.33,135.36) .. controls (230.09,135.34) and (232.35,137.56) .. (232.36,140.33) .. controls (232.38,143.09) and (230.15,145.34) .. (227.39,145.35) -- cycle ;
			\draw  [fill={rgb, 255:red, 65; green, 117; blue, 5 }  ,fill opacity=1 ] (186.39,145.6) .. controls (183.63,145.61) and (181.38,143.39) .. (181.36,140.63) .. controls (181.35,137.87) and (183.57,135.61) .. (186.33,135.6) .. controls (189.1,135.58) and (191.35,137.81) .. (191.36,140.57) .. controls (191.38,143.33) and (189.15,145.58) .. (186.39,145.6) -- cycle ;
			\draw  [fill={rgb, 255:red, 208; green, 2; blue, 27 }  ,fill opacity=1 ] (145.39,145.84) .. controls (142.63,145.86) and (140.38,143.63) .. (140.36,140.87) .. controls (140.35,138.11) and (142.57,135.86) .. (145.33,135.84) .. controls (148.1,135.83) and (150.35,138.05) .. (150.36,140.81) .. controls (150.38,143.57) and (148.16,145.83) .. (145.39,145.84) -- cycle ;
			\draw    (209.65,188.46) -- (267.36,140.12) ;
			\draw    (209.65,188.46) -- (145.36,140.84) ;
			\draw    (209.65,188.46) -- (197.14,162.75) -- (186.36,140.6) ;
			\draw    (209.65,188.46) -- (227.36,140.36) ;
			\draw  [fill={rgb, 255:red, 208; green, 2; blue, 27 }  ,fill opacity=1 ] (145.39,193.84) .. controls (142.63,193.86) and (140.38,191.63) .. (140.36,188.87) .. controls (140.35,186.11) and (142.57,183.86) .. (145.33,183.84) .. controls (148.1,183.83) and (150.35,186.05) .. (150.36,188.81) .. controls (150.38,191.57) and (148.16,193.83) .. (145.39,193.84) -- cycle ;
			\draw [fill={rgb, 255:red, 208; green, 2; blue, 27 }  ,fill opacity=1 ]   (145.36,188.84) -- (145.36,140.84) ;
			\draw    (209.65,188.46) -- (145.36,188.84) ;
			\draw    (209.65,188.46) -- (268.33,189.33) ;
			\draw  [fill={rgb, 255:red, 74; green, 144; blue, 226 }  ,fill opacity=1 ] (268.36,194.33) .. controls (265.6,194.35) and (263.35,192.12) .. (263.33,189.36) .. controls (263.32,186.6) and (265.54,184.35) .. (268.3,184.33) .. controls (271.07,184.32) and (273.32,186.54) .. (273.33,189.3) .. controls (273.35,192.07) and (271.12,194.32) .. (268.36,194.33) -- cycle ;
			\draw  [fill={rgb, 255:red, 248; green, 231; blue, 28 }  ,fill opacity=1 ] (209.36,241.33) .. controls (206.6,241.35) and (204.35,239.12) .. (204.33,236.36) .. controls (204.32,233.6) and (206.54,231.35) .. (209.3,231.33) .. controls (212.07,231.32) and (214.32,233.54) .. (214.33,236.3) .. controls (214.35,239.07) and (212.12,241.32) .. (209.36,241.33) -- cycle ;
			\draw    (209.65,188.46) -- (209.33,236.33) ;
			\draw    (268.33,189.33) -- (209.33,236.33) ;
			\draw    (267.36,140.12) -- (209.33,236.33) ;
			\draw    (227.36,140.36) -- (186.36,140.6) ;
			
			\draw (179.04,112.91) node [anchor=north west][inner sep=0.75pt]  [rotate=-359.12] [align=left] {2};
			\draw (137.27,112.87) node [anchor=north west][inner sep=0.75pt]  [rotate=-359.92] [align=left] {1};
			\draw (219.04,112.91) node [anchor=north west][inner sep=0.75pt]  [rotate=-359.12] [align=left] {3};
			\draw (260.04,112.91) node [anchor=north west][inner sep=0.75pt]  [rotate=-359.12] [align=left] {4};
			\draw (138.04,194.91) node [anchor=north west][inner sep=0.75pt]  [rotate=-359.12] [align=left] {5};
			\draw (209.04,194.91) node [anchor=north west][inner sep=0.75pt]  [rotate=-359.12] [align=left] {6};
			\draw (262.04,194.91) node [anchor=north west][inner sep=0.75pt]  [rotate=-359.12] [align=left] {7};
			\draw (203.04,242.91) node [anchor=north west][inner sep=0.75pt]  [rotate=-359.12] [align=left] {8};

		\end{tikzpicture}
	
	\textbf{G}
	\end{center}
	
	\columnbreak

\begin{center}
		\tikzset{every picture/.style={line width=0.75pt}} 
	
	\begin{tikzpicture}[x=0.75pt,y=0.75pt,yscale=-1,xscale=1]
		
		\draw  [fill={rgb, 255:red, 0; green, 0; blue, 0 }  ,fill opacity=1 ] (454.68,201.56) .. controls (451.92,201.57) and (449.66,199.35) .. (449.65,196.59) .. controls (449.63,193.82) and (451.86,191.57) .. (454.62,191.56) .. controls (457.38,191.54) and (459.63,193.76) .. (459.65,196.53) .. controls (459.66,199.29) and (457.44,201.54) .. (454.68,201.56) -- cycle ;
		\draw  [fill={rgb, 255:red, 74; green, 144; blue, 226 }  ,fill opacity=1 ] (512.39,153.21) .. controls (509.63,153.23) and (507.38,151) .. (507.36,148.24) .. controls (507.35,145.48) and (509.57,143.23) .. (512.33,143.21) .. controls (515.09,143.2) and (517.35,145.42) .. (517.36,148.18) .. controls (517.38,150.94) and (515.15,153.2) .. (512.39,153.21) -- cycle ;
		\draw  [fill={rgb, 255:red, 65; green, 117; blue, 5 }  ,fill opacity=1 ] (431.39,153.69) .. controls (428.63,153.71) and (426.38,151.48) .. (426.36,148.72) .. controls (426.35,145.96) and (428.57,143.71) .. (431.33,143.69) .. controls (434.1,143.68) and (436.35,145.9) .. (436.36,148.66) .. controls (436.38,151.43) and (434.15,153.68) .. (431.39,153.69) -- cycle ;
		\draw  [fill={rgb, 255:red, 208; green, 2; blue, 27 }  ,fill opacity=1 ] (390.39,153.94) .. controls (387.63,153.95) and (385.38,151.73) .. (385.36,148.97) .. controls (385.35,146.2) and (387.57,143.95) .. (390.33,143.94) .. controls (393.1,143.92) and (395.35,146.15) .. (395.36,148.91) .. controls (395.38,151.67) and (393.16,153.92) .. (390.39,153.94) -- cycle ;
		\draw    (454.65,196.56) -- (512.36,148.21) ;
		\draw    (454.65,196.56) -- (390.36,148.94) ;
		\draw    (454.65,196.56) -- (442.14,170.84) -- (431.36,148.69) ;
		\draw  [fill={rgb, 255:red, 248; green, 231; blue, 28 }  ,fill opacity=1 ] (454.36,249.43) .. controls (451.6,249.44) and (449.35,247.22) .. (449.33,244.46) .. controls (449.32,241.7) and (451.54,239.44) .. (454.3,239.43) .. controls (457.07,239.41) and (459.32,241.64) .. (459.33,244.4) .. controls (459.35,247.16) and (457.12,249.41) .. (454.36,249.43) -- cycle ;
		\draw    (454.65,196.56) -- (454.33,244.43) ;
		\draw    (512.36,148.21) -- (454.33,244.43) ;
		
		\draw (424.04,121) node [anchor=north west][inner sep=0.75pt]  [rotate=-359.12] [align=left] {2};
		\draw (382.27,120.96) node [anchor=north west][inner sep=0.75pt]  [rotate=-359.92] [align=left] {1};
		\draw (505.04,121) node [anchor=north west][inner sep=0.75pt]  [rotate=-359.12] [align=left] {4};
		\draw (454.04,203) node [anchor=north west][inner sep=0.75pt]  [rotate=-359.12] [align=left] {6};
		\draw (448.04,251) node [anchor=north west][inner sep=0.75pt]  [rotate=-359.12] [align=left] {8};

	\end{tikzpicture}

	\textbf{H}
\end{center}

\end{multicols}


An example of a graph $G$ and its reduced graph $H$ is given in the above picture. Also, the vertices of $G$ which belong to the same twin classes are colored by the same colour. We observe that graph $G$ has eight vertices whereas graph $H$ has only five vertices.  
	
\end{rem}

On the other hand, topological indices are parameters associated with graphs that are proven to have many applications in different areas such as mathematical chemistry \cite{w}. A list of many topological indices and their application can be seen in \cite{aja2019, survey1}. Among various topological indices, the Wiener index is classical and an important topological index. 

Steiner distance is a well-studied notion in the literature which is also shown to have applications in different fields such as networking \cite{d20, mk09,oo99, fh00, pho97, lo08, woh94}. Steiner distance is studied from a theoretical point of view as well. For example, recently the spectral aspects of the Steiner distance matrix are explored in \cite{abs22}.

In this paper, using the above interpretation of the graph $G$ as the generalized composition of subgraphs induced by the twin classes, we will derive a formula for the $m$-Steiner Wiener index of the graph $G$. In particular, we will show that calculating this parameter for a graph $G$ can be reduced to calculating the same for a smaller graph $H$.






We prove the following formula for the $m$-Steiner Wiener index of an arbitrary graph $G$ (in terms of the cardinalities of the twin classes). 

\begin{theom}
	Let $G$ be an arbitrary finite connected graph with the associated twin classes $C_1,\dots,C_k$. Then the  $m$-Steiner Wiener index of $G$ is given by \begin{equation}\label{main}
		\small	SW_m(G) = \sum_{\substack{1 \le i \le k \\ G_i \text{ is complete}}}(m-1) \binom{n_i}{m} + \sum_{\substack{1 \le i \le k \\ G_i \text{ is empty}}}m \binom{n_i}{m}+\sum_{\substack{S \in \mathcal P_m(V) \\ |I(S)|=l>1}} \Big(d_H^l(S_r)+ \sum\limits_{c_i \in S_r}(t_i - 1))\Big).
	\end{equation} where $|C_i| = n_i$. 

\end{theom}

	If we try to calculate $SW_m(G)$ from its definition then we need to calculate the $m$ Steiner distance of all the $m$ element subsets in the graph $G$, But using our formula,  we have reduced this calculation to calculate the $l$-Steiner distances (always $l \le m$) of the smaller subsets in the graph $H$. In general, the graph $H$ of representatives will be much smaller than the graph $G$, and hence calculating these distances in this smaller graph $H$ is much simpler. This is illustrated in Example \ref{ezn}. This explains the usefulness of our formula given in the above theorem.

In the subsequent sections, as applications of the above theorem, the $m$-Steiner Wiener index is calculated for various important classes of graphs from the literature  including

1. Power graphs associated with groups

2. Zero divisor graphs and the ideal-based zero divisor graphs associated with commutative rings with unity.

3. Comaximal ideal graphs associated with commutative rings with unity.

This explains the applicability of our formula to various cases. Many examples are provided throughout the paper explaining our results.

\begin{rem}
	The Wiener index is one of the important examples of topological indices that we have calculated in this paper. Using our approach in this paper, we can find a formula for any topological indices in terms of cardinalities of twin classes. Since finding such a formula for other topological indices requires the same work as in this paper, we deal only with the Wiener index and Steiner Wiener index in this paper. 
\end{rem}

\section{A formula for the Wiener index of graphs using twin classes}
In what follows, we use the graph-theoretic notions from \cite{dbw}. In this section, we define the notion of twin classes of a graph $G$. Then we will show that the graph $G$ can be visualized as a generalized composition of subgraphs induced by the twin classes \cite{bst93, tfz08}. Using this result, we derive a formula for the $m$- Steiner Wiener index of the graph $G$.

\subsection{Generalized composition of graphs}

The generalized composition of graphs, denoted by $H[G_1,G_2,\dots,G_k]$, is defined as follows \cite{sch}.
\begin{defn}\label{gjoin}\cite[Equation (48) below]{sch}
	Let $H$ be a graph with vertex set $V = \{1,\dots,k\}$ and let $G_1,\dots,G_k$ be a collection of graphs with the respective vertex sets $V_i = \{v_i^{1},\dots,v_i^{n_i}\}$ for $(1 \le i \le k)$. Then their generalized composition $G = H[G_1,\dots,G_k]$ has vertex set $V_1 \sqcup \cdots \sqcup V_k$ and two vertices $v_i^p$ and $v_j^q$ of $H[G_1,\dots,G_k]$ are adjacent if one the following conditions is satisfied:
	\begin{enumerate}
		\item $i=j$, and $v_i^p$ and $v_j^q$ are adjacent vertices in $G_i$.
		\item $i \ne j$ and $i$ and $j$ are adjacent in $H$.
	\end{enumerate}
\end{defn}
	The graph $H$ is said to be the base graph of $G$ and the graphs $G_i$'s are called the factors of $G$. In the next subsection, we will prove that any arbitrary graph $G$ is a generalized composition of graphs induced by the twin classes. 


\subsection{Representing arbitrary graphs as generalized composition using twin classes}\label{zero as join}

The results given in this subsection are well-known \cite{bst93, tfz08}. For the exposition, we have included them here. 

Let $R$ be a finite commutative ring with unity. In \cite[Section 3]{al99}, when $R = \mathbb Z_n$ the ring of integers modulo $n$, the zero-divisor graph of $R$ is realized as a generalized composition of suitable choices of graphs. Extending this result, in \cite[Section 2.3]{sga}, for an arbitrary ring $R$, the zero-divisor graph $\Gamma(R)$ of $R$ is realized as a generalized composition of suitable choices of graphs. In the same paper, using this realization, a formula for the Wiener index of $\Gamma(R)$ is proved.  In this subsection,  we generalize this realization to an arbitrary graph $G$ by representing it as a generalized composition using the concept of twin classes. Further, by using this realization we prove a formula for the Wiener index of the graph $G$ in terms of the sizes of the twin classes of $G$. Indeed, we prove a formula for the more general $m$-Steiner Wiener index of the graph $G$. 


Let $G$ be an arbitrary finite graph with vertex set $V(G)$. For a vertex $v$ of $G$ the open neighbourhood of $v$, denoted by $N_G(v)$, is defined to be the set of vertices of $G$ which are adjacent to $v$. Similarly the closed neighbourhood of $v$ is defined to be $N_G[v]: = N_G(v) \cup \{v\}$. Two vertices $u$ and $v$ of $G$ are said to be \emph{true twin vertices} if $N_G[u] = N_G[v]$ and \emph{false twin vertices} if $N_G(u) = N_G(v)$. 
In general, the vertices $u$ and $v$ are said to be twin vertices if they are either true twins or false twins.

We define an equivalence relation $\sim^G$ on $V(G)$ as follows. For $u,v \in V(G)$, define $u \sim^G v$ if, and only if, $N_G(u) \backslash \{v\} = N_G(v) \backslash \{u\}$. In other words, two vertices of $G$ are related by $\sim^G$ if they are twin vertices in $G$. Let $C_1, \dots, C_k$ be the equivalence classes of this relation (with respective representatives $c_1,\dots,c_k$) which we call the twin classes of the graph $G$.

The following lemma is proved in \cite[Lemma 3.1]{tfz08} which illustrates the structure of the graph $G$ in terms of the twin classes. 

\begin{lem}\label{a c or t}
	\begin{enumerate}
		\item The induced subgraph $G(C_i)$ of the graph $G$ is either complete or totally disconnected. 
		\item If $i \ne j$, a vertex $a$ of $C_i$ is adjacent to either all or none of the vertices of $C_j$.
	\end{enumerate}
\end{lem}


Let $H$ be the subgraph induced by the vertices $\{c_1,\dots,c_k\}$ in $G$. $H$ is called the reduced graph of the graph $G$. This smaller graph $H$ reflects many structural properties of the bigger graph $G$. In particular, we have the following result.


\begin{lem}
	Let $G$ be a graph with at least two twin classes. Then the graph $G$ is connected if, and only if, the graph $H$ is connected.
\end{lem}
\begin{pf}
	The proof follows from the more general statement given below.
\end{pf}

The following proposition says that $G$ is a generalized composition of certain complete graphs and totally disconnected graphs.

\begin{prop}\label{join}
	Let $G$ be an arbitrary graph possibly with loops. Let the twin classes of the graph $G$ be $C_1,\dots,C_k$. Define $G_i := G[C_i]$ and $H := G[\{c_1,\dots,c_k\}]$. Then $G = H[G_1,\dots,G_k]$.
\end{prop}
\begin{pf}
	The proof follows from Lemma \ref{a c or t}.
\end{pf}


\subsection{The  Steiner Wiener index of generalized composition of graphs} \label{ks}
In \cite[Section 2.4]{sga}, a formula for the Wiener index of the generalized composition of graphs is proved. In this subsection, we extend this formula to $m$-Steiner Wiener index using twin classes.
\begin{defn}\label{sw2}
	Let $G$ be a connected graph with vertex set $V$ and let $S \in \mathcal P_m(V)$ the set of all $m$-element subsets of the set $V$. The $m$-Steiner distance of the set $S$, denoted by $d_G^m(S)$, is defined to be the number of edges in the smallest subtree of $G$ containing $S$. The \emph{ m-Steiner Wiener index} of the graph $G$, denoted by $SW_m(G)$,  is defined to be $$SW_m(G) = \sum_{\substack{S \subseteq V \\ |S|=m}} d_G^m(S).$$
	\end{defn}

Let $G$ be a connected graph with twin classes $C_1,\dots,C_k$ with respective representatives $c_1,\dots,c_k$ and $G = H[G_1,\dots,G_k]$ where $H = \{c_1,\dots,c_k\}$ and the factor graphs $G_i = G[C_i]$'s are either complete or empty (graph with vertices but no edges). Let $V(G_i) = \{v_i^{1},\dots,v_i^{n_i}\}$ be the vertex set of the graph $G_i$ and $V(G) = \sqcup_{i=1}^kV(G_i)$. Given these, we define the following notions.


\begin{defn}
  Let $S \in \mathcal P_m(V(G))$. The intersection number of $S$ is defined to be the cardinality of the set $I(S) := \{1 \le i \le k : S \cap C_i \ne \phi\}$. We also define the set $S_r = \{c_i : S \cap C_i \ne \phi\}$ which consists of one representative of each twin class which intersects $S$. We also define $t_i := |S \cap C_i|$ for $1 \le i \le k$.
\end{defn}

The following lemma explains the Steiner distance of certain sets lying in both the graphs $G$ and $H$.
\begin{lem}\label{dist H and G}
	Let $S$ be a set of $l$ vertices of $G$ such that the set $S \cap C_i$ has at most one element for each $1 \le i \le k$. Then $S = S_r$ (by changing the representatives $c_i$ if necessary) and $d_H^l(S) = d_G^l(S)$.
\end{lem}
\begin{pf}
	The proof follows from the definition of the graph $G$. 
\end{pf}

 The following lemma is an extension of \cite[Lemma 3]{sga} which explains the $m$-Steiner distance of any set $S$ of $m$ vertices of $G$.

\begin{lem} \label{general distance} Let $S \in \mathcal P_m(V(G))$.  Then the $m$-Steiner distance of the set $S$ in the graph $G$ is given by 
	\begin{equation}\label{dist}d_G^m(S) = \begin{cases}
			m-1 &\text{ if } I(S)=\{i\} \text{ and }G_i \text{ is complete}\\
			m &\text{ if } I(S)=\{i\} \text{ and }G_i \text{ is totally disconnected}\\
			d_H^l(S_r)+ \sum\limits_{c_i \in S_r}(t_i - 1) &\text{ if } l:=|I(S)|>1. 
	\end{cases}\end{equation} 
\end{lem}
\begin{pf}
Fix $S \in \mathcal P_m(V(G))$ and let $S_r = \{c_{i_1},\dots,c_{i_l}\}$. Without loss of generality, by changing the representatives of the classes $C_i$ if required, we can assume that $S_r \subseteq S$. 

Suppose $S_r = \{c_i\}$ and $G_i$ is complete then the proof is immediate. 

Suppose $S_r = \{c_i\}$ and $G_i$ is totally disconnected. Let $c_j \in V(H)$ be such that $c_j$ is adjacent to $c_i$.  Now, every element of $S$ ($=S \cap C_i$) is adjacent to $v_j^q$ where $1 \le q \le n_j$. Therefore, the vertices $S \cup \{v_j^q\}$ form a star graph on $m+1$ vertices and the result follows.

Suppose $|S_r| = l>1$. Let $T$ be the smallest subtree of $G$ containing $S_r$. Then, by the above lemma, $d_H^l(S_r) = d_G^l(S_r) = h$ (say) where $h$ is the number of edges in $T$. Let $c_i \in S_r$. Then there are $t_i= |S \cap C_i|$ many vertices of $C_i$ are in $S$. If $G_i$ is a complete graph then we can consider the star subgraph (not necessarily induced) of $G_i$ formed by  these $t_i$ many vertices whose middle vertex is $c_i$. i.e., they can be joined by $t_i-1$ edges. Also, we can clue this star graph to the tree $T$ at the common vertex $c_i$. If $G_i$ is totally disconnected then consider the vertex $c$ of $T$ which is adjacent to $c_i$. Since $c_i$ is adjacent to $c$, every vertex in $S \cap C_i$ is adjacent to $c$. This way we get a star subgraph of $G$ formed by the $t_i$ many vertices $((S \cap C_i) \backslash \{c_i\} \cup \{c\})$ whose middle vertex is $c$. i.e., they can be joined again by $(t_i-1)$ edges. Also, we can clue this star graph to the tree $T$ at the common vertex $c$. Note that, the vertices $c$ and $c_i$ are already adjacent in $T$ and their distance is accounted for in $d_H^l(S_r)$.  

By doing this process for each $c_i \in S_r$, the tree $T$ is extended to a tree $T^{'}$ of $G$ which by construction the smallest subtree containing the set $S$ with the required number of edges. This completes the proof.
\end{pf}

From the above lemma, we get the following result.

\begin{thm}\label{mainthm} 
	Let $G$ be an arbitrary connected finite graph. With the notations defined above,  the  $m$-Steiner Wiener index of $G$ is given by \begin{equation}\label{main}
	\small	SW_m(G) = \sum_{\substack{1 \le i \le k \\ G_i \text{ is complete}}}(m-1) \binom{n_i}{m} + \sum_{\substack{1 \le i \le k \\ G_i \text{ is empty}}}m \binom{n_i}{m}+\sum_{\substack{S \in \mathcal P_m(V) \\ |I(S)|=l>1}} \Big(d_H^l(S_r)+ \sum\limits_{c_i \in S_r}(t_i - 1)\Big).
	\end{equation} where $|C_i| = n_i$.

\end{thm}
\begin{pf}
	The proof follows from the above discussions.
\end{pf}

\begin{cor}
	Let $S=\{v_i^{p},v_j^q\} \in \mathcal P_2(V)$ such that $|I(S)|=2$.  Then we have $i \ne j$, $t_i = t_j=1$ and $d_G^2(S) = d_H^2(S)=d_H(i,j)$ by Lemma \ref{dist H and G}. 
	
	Now, if we assume that $m=2$ in the above theorem, we get the Wiener index of the graph $G$ is equal to
	\begin{align}
		\small	SW_2(G) &= \sum_{\substack{1 \le i \le k \\ G_i \text{ is complete}}} \binom{n_i}{2} + \sum_{\substack{1 \le i \le k \\ G_i \text{ is empty}}}2 \binom{n_i}{2}+\sum_{\substack{S \in \mathcal P_2(V) \\ |I(S)|=2}} \Big(d_H^l(S_r)+ \sum\limits_{c_i \in S_r}(t_i - 1)\Big) \\
		&= \sum_{\substack{1 \le i \le k \\ G_i \text{ is complete}}} \binom{n_i}{2} + \sum_{\substack{1 \le i \le k \\ G_i \text{ is empty}}}2 \binom{n_i}{2}+\sum_{\substack{S \in \mathcal P_2(V) \\ |I(S)|=2}} d_H(i,j). \\
		&= \sum_{\substack{1 \le i \le k \\ G_i \text{ is complete}}} \binom{n_i}{2} + \sum_{\substack{1 \le i \le k \\ G_i \text{ is empty}}}2 \binom{n_i}{2}+\sum_{1 \le i < j \le k} n_i n_j d_H(i,j). \\
	\end{align} 
 which is the formula for the Wiener index of an arbitrary graph $G$ in terms of the cardinalities of its twin classes. 
 
\end{cor}

\begin{rem}\label{remark}
	From the above theorem, if we know the cardinalities of the twin classes $C_i$ and the distances between vertices in $H$ then the $m$-Steiner Wiener index of the arbitrary graph $G$ can be explicitly determined. 
	
\end{rem} 

In the following sections, as applications of Theorem \ref{mainthm}, we calculate the $m$-Steiner Wiener index of the graph $G$ when

1. $G$ is the $K_n$-join of some graphs.

2. $G$ is the power graph of an arbitrary finite group $\Gamma$.

3. $G$ is the zero-divisor graph of an arbitrary finite commutative ring with unity. Indeed, we work with a more general ideal-based zero divisor graph of a finite commutative ring with unity and

4. Comaximal ideal graph of a finite commutative ring with unity.

\section{The Steiner Wiener index of $K_n$-join graphs}
In this section, using our twin class and the $H$-join approach discussed in the previous section, we study the $m$-Steiner Wiener index of an infinite class of graphs called the completely joined graphs. For this class of graphs, we give an upper bound on its $m$-Steiner Wiener index. As a corollary, we explicitly calculate the $m$-Steiner Wiener index of the  complete multipartite graphs.
\begin{defn}
	Any graph $G = K_p[G_1,\dots,G_p]$ where $G_i$s are arbitrary graphs is said to be a completely joined graph.
	
\end{defn}
\begin{example}
	\begin{enumerate}
		\item Complete bipartite ($K_2[\overline{K_{n_1}},\overline{K_{n_2}}]$) and Complete multipartite graphs ($K_p[\overline{K_{n_1}},\dots,\overline{K_{n_p}}]$) are examples of completely joined graphs. 
		\item Any graph $G$ which is a join of two graphs $G_1$ and $G_2$, i.e., $G = G_1 \vee G_2 = K_2[G_1,G_2]$, is a completely joined graph. In particular, wheel graphs ($K_2[C_n,K_1]$ for $n \ge 3$) and the special class of split graphs $K_2[K_n,\overline{K_m}]$ are completely joined graphs.
		\item Let $G$ be a graph with vertex set $V(G)$.  As in Section \ref{zero as join} write $G = H[G_1,\dots,G_p]$ where $G_i = G[C_i]$ and $H = G[\{c_1,\dots,c_p\}]$. Given this, graph $G$ is completely joined graph if the reduced graph $H$ is complete. 
	\end{enumerate}
\end{example}

Let $G=K_p[G_1,\dots,G_p]$ (where $G_i$s are arbitrary graphs) be a completely joined graph. In the following lemma, we prove an upper bound for the $m$-Steiner distance of an $m$-element subset of $V(G)$.
\begin{lem}
	Let $S$ be a subset of $V(G)$ such that $|S| = m$. Then 
	
	$$d_G^m(S) = \begin{cases}
		\text{ atmost }m \text{    if } |I(S)| = 1 \\
		m-1 \text{ if } |I(S)|>1 
	\end{cases}.$$
where $I(S) = \{1 \le i \le p : S \cap V(G_i) \ne \phi\}$.
\end{lem}
\begin{pf}
	Let $S$ be a subset of $V(G)$ such that $|S| = m$. 
	
	Assume $|I(S)|>1$. Fix $s \in S$. Then $s \in V(G_i)$ for some $i$ and this $s$ is adjacent to $v$ for all $v \in S \cap V(G_j)$ for $j \ne i$.  This way we can form a star subgraph $T$ of $G$ of size $|S|-t_i+1$. Consider the graph $G_j$ such that $V(G_j) \cap S \ne \phi$ and $j \ne i$. Now, $|S \cap V(G_j)| = t_j$ and all these $t_j$ many vertices are adjacent to $s$. Fix a vertex $q$ in $S \cap V(G_j)$. Then every vertex in $(S \cap V(G_i)) \backslash \{s\}$ is adjacent to $q$.
	Therefore the vertices in  $(S \cap V(G_i)) \backslash \{s\}$ can be added to the star subgraph $T$ to form a bigger tree $T^{'}$ whose vertex set is precisely $S$. By the construction of $T^{'}$, it has $m-1$ edges.
	
	Assume that $|I(S)|=1$. Then $S \subseteq V(G_i)$ for some $i$. Fix a vertex $q$ in the graph $G_j$ where $j \ne i$. Then $q$ is adjacent to every element of $S$ and hence we can form a subtree $T$ of $G$ on $m+1$ vertices. Note that this tree $T$ need not be the smallest one.  Therefore $d_G^m(S) \le m$ in this case. This completes the proof of the lemma.
\end{pf}

In the following theorem, we have given an upper bound on the $m$-Steiner Wiener index of the completely joined graphs. The proof of this theorem follows from the above lemma.
\begin{thm}
	Let $G=K_p[G_1,\dots,G_p]$ where $G_i$s are arbitrary graphs be a completely joined graph. Then $$SW_m(G) \le \binom{n}{m}(m)$$ where $n = \sum_{i=1}^{p}|V(G_i)|$. 
\end{thm}

\begin{rem}
	We observe that the above bound on the $m$-Steiner Wiener index of $G$ depends only on the number of vertices in $G$ and the value of $m$. In particular, it is independent of the structure of the graphs $G_i$ for $1 \le i \le p$.
\end{rem}

\begin{example}
	Consider the wheel graph $W_n := K_2[C_n,K_1]$ for $n \ge 3$.  Let $2 \le m \le n+1$, then $SW_m(W_n) \le \binom{n+1}{m} (m)$. In particular, the Wiener index of the graph $W_n$ can be atmost $2 \binom{n+1}{2} = (n+1)(n)$. 
\end{example}

In the following corollary, we will explicitly find the $m$-Steiner Wiener index of the complete multipartite graph $K_{m_1,\dots,m_p}$. The following lemma will be useful for the same. 

\begin{lem}
	Consider the complete multipartite graph $G = K_{m_1,\dots,m_p}$. Let $S$ be a size $m$ subset of the vertex set of $G$. Then $d_G^m(S) = \begin{cases}
		m \text{    if } |I(S)| = 1 \\
		m-1 \text{ if } |I(S)|>1 
	\end{cases}$.
\end{lem}
\begin{pf}
	Let $G = K_{m_1,\dots,m_p}$ be the complete multipartite graph and $S$ be a size $m$ subset of the vertex set of $G$. Note that $G = K_p[G_1,\dots,G_p]$ where $G_i$ is the graph $\overline{K_{n_i}}$. We assume that $t_i = |S \cap V(G_i)|$ and hence $|S| = \sum_{i=1}^{p}t_i$. Assume that $|I(S)|=1$. Then $S \subseteq V(G_i)$ for some $i$. Fix a vertex $q$ in the graph $G_j$ where $j \ne i$. Then $q$ is adjacent to every element of $S$ and hence we can form a subtree $T$ of $G$ on $m+1$ vertices. Since $G_i$ is totally disconnected and $S \subseteq V(G_i)$ there cannot be a subtree of $G$ with $m$ vertices. Therefore $d_G^m(S) = m$ in this case.

	Assume $|I(S)|>1$. Fix $s \in S$. Then $s \in V(G_i)$ for some $i$ and this $s$ is adjacent to $v$ for all $v \in S \cap V(G_j)$ for $j \ne i$.  This way we can form a star subgraph $T$ of $G$ of size $|S|-t_i+1$. Consider the graph $G_j$ such that $V(G_j) \cap S \ne \phi$ and $j \ne i$. Now, $|S \cap V(G_j)| = t_j$ and all these $t_j$ many vertices are adjacent to $s$. Fix a vertex $q$ in $S \cap V(G_j)$. Then every vertex in $(S \cap V(G_i)) \backslash \{s\}$ is adjacent to $q$.
	Therefore the vertices in  $(S \cap V(G_i)) \backslash \{s\}$ can be added to the star subgraph $T$ to form a bigger tree $T^{'}$ whose vertex set is precisely $S$. By the construction of $T^{'}$, it has $m-1$ edges. This completes the proof.
\end{pf}
\begin{cor}
	From the above lemma, the $m$-Steiner Wiener index of the complete multi partite graph $G=K_{n_1,\dots,n_p}$ is equal to $$SW_m(G) = (\binom{n}{m}-\sum_{i=1}^{p}\binom{n_i}{m}) (m-1) +  \sum_{i=1}^{p}\binom{n_i}{m} (m)  = \binom{n}{m}(m-1)+\sum_{i=1}^{p}\binom{n_i}{m}$$
where $n = \sum_{i=1}^{p}m_i$. In particular, when $m=2$, the Wiener index of the complete multipartite graph $G=K_{m_1,\dots,m_p}$ is equal to $\binom{n}{2} + \sum_{i=1}^{p}\binom{n_i}{2}$ where $n = \sum_{i=1}^{p}m_i$.
\end{cor}

For example, consider the graph $G = K_{3,3,3}$ and $m=5$. Then 
$SW_5(G) = \binom{9}{5}(4) = 504$.

 \section{The Steiner Wiener index of power graphs of groups}\label{power graph}

 Let $\Gamma$ be an arbitrary finite group. Let $G$ be the power graph of the group $\Gamma$.
 In this section, as an application of Theorem \ref{mainthm}, we derive a formula for the $m$-Steiner Wiener index of the graph $G$. As further applications, we have explicitly calculated the $m$-Steiner Wiener index of the power graph of various classes of finite groups including
 
 1. group of integers modulo $n$, 
 
 2. dihedral group and 
 
 3. the quaternion group is discussed. 
 
 We start with the definition of the power graph of a finite group.
 \begin{defn}
 	Let $\Gamma$ be a finite group. The power graph $G$ of the group $\Gamma$ has vertex set $\Gamma$ and two group elements $a$ and $b$ are adjacent in $G$ if either $a$ is a power of $b$ or $b$ is a power of $a$.
 \end{defn}

Let $\Gamma$ be a finite group and $G$ be the associated power graph.
We need to understand the twin classes of $G$ to apply Theorem \ref{mainthm} to the graph $G$. For a group element $a$, we define $b = a^{\frac{1}{k}}$ for some $k \in \mathbb N$ if $b^k = a$ and we define 

$$Pow(a) := \{b \in \Gamma : b = a^n \text{ for some } n \in \mathbb N \cup \{\frac{1}{k}: k \in \mathbb N\}\}.$$ 


Note that, the set $Pow(a)$ need not be equal to the subgroup generated by $a$, and the elements of $Pow(a)$ are the elements of $\Gamma$ that are adjacent to $a$ in $G$. Define an equivalence relation $\sim$ on $V(G)$ by $a \sim b$ if $Pow(a)\backslash \{b\} = Pow(b)\backslash \{a\}$.  Let $C_1, \ldots,C_t$ be the twin classes of this relation with respective representatives $c_1,\ldots,c_t$. 

We denote by $H(G)$ the subgraph of $G$ induced by the representatives $c_1, \dots,   c_t$. Then by Proposition \ref{join} we have




\begin{prop}\label{p join}
	Let $G$ be the power graph of a finite group $\Gamma$. Let $G[C_i]$ be the subgraph induced by the twin classes $C_i$ in $G$. Then $G = H[G[C_1],\dots,G[C_t]]$ 
\end{prop}



Now, we can use Theorem  \ref{mainthm} to calculate the $m$-Steiner Wiener index of $G$.

\begin{thm}\label{ppp}
	Let $G$ be the power graph of a finite group $\Gamma$. Then the  $m$-Steiner Wiener index of $G$ is given by

	{\small \begin{equation}
			SW_m(G)) = \sum_{\substack{1 \le i \le t \\ \mathfrak{g}(C_i) \text{ is complete}}}(m-1) \binom{n_i}{m} + \sum_{\substack{S \in \mathcal P_m(V) \\ |I(S)|=l>1}} \Big(d_H^l(S_r)+ \sum\limits_{c_i \in S_r}(t_i - 1)\Big).
	\end{equation}} where $|C_i| = n_i$.
	
\end{thm}

\begin{cor}\label{pw}
	Let $G$ be the power graph of a finite group $\Gamma$. Then the  Wiener index of $G$ is given by \begin{equation}
		W(G) = \sum_{\substack{1 \le i \le k \\ \mathfrak{g}(C_i) \text{ is complete}}}\binom{n_i}{2} +\sum_{1 \le i<j \le k}n_in_jd_H(i,j).
	\end{equation} where $|C_i| = n_i$.  
\end{cor}

In the following subsections, using the above theorem, we calculate the $m$-Steiner Wiener index of power graphs explicitly for various classes of finite groups.


Let $\Gamma$ be the group of integers modulo $n$.  Let $G$ be the power graph of $\Gamma$. First, we write the graph $G$ as a $H$-join of some suitable choices of induced subgraphs. Let $d$ be a  divisor of $n$ such that $d \ne n$. The following sets $V_d$ are defined in \cite[Section 3]{al99}.  $$V_d := \begin{cases} \{x \in \Gamma : \gcd(x,n) = d \} &\text{ if } d>1 \\ \{x \in \Gamma : \gcd(x,n) = d \}\cup \{0\} &\text{ if } d=1 \end{cases}.$$

We make the following observations on the set $V_d$.

\begin{enumerate}
	\item The number of elements in the set $V_{d}$ is equal to $\phi(\frac{n}{d})$ if $d>1$ and $\phi(n)+1$ if $d=1$.
	\item $a \in V_d$ implies that order of $a$ in the group $\Gamma$ is equal to $\frac{n}{d}$.
	\item For a divisor $d$ of $n$, the subgraph $G_d$ which is defined to be the subgraph induced by the set $V_d$ in $G$ is a complete graph.
	\item The graph $H$ defined by $V(H)$ is the set of all divisors of $n$ and $E(H) = \{d_id_j : d_i|d_j \text{ or } d_j|d_i\}$ is an induced subgraph of $G$.
	\item Suppose $d_i | d_j$ then $o(x) | o(y)$ for every $x \in V_{d_i}$ and $y \in V_{d_j}$. In particular, every element of $V_{d_i}$ is adjacent to every element of $V_{d_j}$ in the graph $G$.
	\item The twin classes of the graph $G$ are precisely the sets $V_d$ defined above.
	
\end{enumerate}

Given the above observations, we have the following proposition
\begin{prop}\label{pzn}
	The power graph $G$ can be written as $H[G_d: \text{ $d < n$ is a divisor of n }]$ where the graphs $H$ and $G_d = G[V_d]$ are defined above.  In particular $$G = H[K_{\phi(n)+1}, K_{\phi(d_1)}, K_{\phi(d_2)}\dots, K_{\phi(d_t)}].$$
	where $1, n, d_1, \dots, d_t$ are the divisors of $n$.
\end{prop}

Now, we can use Theorem \ref{ppp} to calculate the $m$-Steiner Wiener index of $G$ as explained in the next example.

\begin{example} \label{ezn}
	Consider the group $\mathbb{Z}_6=\{\overline{0},\overline{1}, \overline{2}, \overline{3}, \overline{4}, \overline{5} \}$ whose power graph is $G$. We will calculate the  $m$-Steiner Wiener index of $G$ using Theorem \ref{ppp}.  Clearly, $V_1=\{\overline{0}, \overline{1}, \overline{5}\}$,
	$V_2=\{\overline{2},\overline{4}\}$ and $V_3=\{\overline{3}\}$.  Therefore, by Theorem \ref{ppp}, we have
	
		 \begin{equation}
			SW_3(G)) = \sum_{d \text{ is a divisor of 6}}2 \binom{n_i}{3} + \sum_{\substack{S \in \mathcal P_3(V) \\ |I(S)|=l>1}} \Big(d_H^l(S_r)+ \sum\limits_{c_i \in S_r}(t_i - 1)\Big).
	\end{equation} where $|V_i| = n_i$.
	
	\begin{table}[]
		\begin{tabular}{|c|c|l|c|c|}
			\hline
			S.No & 3-subsets & $S_r$                           & $d_H^l(S_r)$ & $\sum\limits_{c_i \in S_r}(t_i - 1)$ \\ \hline
			1    & \{0,1,2\} & \{6,1,2\}                      & 2                             & 0                                                                                     \\ \hline
			2    & \{0,1,3\} & \{6,1,3\}                      & 2                             & 0                                                                                     \\ \hline
			3    & \{0,1,4\} & \{6,1,2\}                      & 2                             & 0                                                                                     \\ \hline
			4    & \{0,1,5\} & \{6,1\}                        & 1                             & 1                                                                                     \\ \hline
			5    & \{0,2,3\} & \{6,2,3\}                      & 2                             & 0                                                                                     \\ \hline
			6    & \{0,2,4\} & \{6,2\}                        & 1                             & 1                                                                                     \\ \hline
			7    & \{0,2,5\} & \{6,1,2\}                      & 2                             & 0                                                                                     \\ \hline
			8    & \{0,3,4\} & \{6,3,2\}                      & 2                             & 0                                                                                     \\ \hline
			9    & \{0,3,5\} & \{6,1,3\}                      & 2                             & 0                                                                                     \\ \hline
			10   & \{0,4,5\} & \{6,1,2\}                      & 2                             & 0                                                                                     \\ \hline
			11   & \{1,2,3\} & \{1,2,3\}                      & 2                             & 0                                                                                     \\ \hline
			12   & \{1,2,4\} & \{1,2\}                        & 1                             & 1                                                                                     \\ \hline
			13   & \{1,2,5\} & \{1,2\}                        & 1                             & 1                                                                                     \\ \hline
			14   & \{1,3,4\} & \multicolumn{1}{c|}{\{1,2,3\}} & 2                             & 0                                                                                     \\ \hline
			15   & \{1,3,5\} & \multicolumn{1}{c|}{\{1,3\}}   & 1                             & 1                                                                                     \\ \hline
			16   & \{1,4,5\} & \multicolumn{1}{c|}{\{1,2\}}   & 1                             & 1                                                                                     \\ \hline
			17   & \{2,3,4\} & \multicolumn{1}{c|}{\{2,3\}}   & 2                             & 1                                                                                     \\ \hline
			18   & \{2,3,5\} & \multicolumn{1}{c|}{\{1,2,3\}} & 2                             & 0                                                                                     \\ \hline
			19   & \{2,4,5\} & \multicolumn{1}{c|}{\{1,2\}}   & 1                             & 1                                                                                     \\ \hline
			20   & \{3,4,5\} & \multicolumn{1}{c|}{\{1,2,3\}} & 2                             & 0                                                                                     \\ \hline
		\end{tabular}
	\caption{}
	\label{tab}
	\end{table}
	
	From Table \ref{tab}, we get 
	$SW_3(G) = 41$. 
	
	
\end{example}


\begin{rem}
	If we try to calculate $SW_3(G)$ from its definition then we need to calculate the $3$ distance of all the sets listed in the second column of the above table in the graph $G$. But using our formula,  we have reduced this calculation to calculate the $3$-distance of the same set of sets in the graph $H$. In general, the graph $H$ of representatives will be much smaller than the graph $G$, and hence calculating the distance in this smaller graph $H$ is much simpler. Given the simplicity of calculating the entries in the last column of the above table, the above example explains the usefulness of our formula given in Theorem \ref{ppp}.  
\end{rem}

	

\begin{example}
	The dihedral group of order $2n$ ($n \ge 3$) is given by
	$$D_{2n} =\left\langle r,s: r^n=1, s^2=1, r^is=sr^{-i}, 0\leq i\leq n-1\right\rangle.$$ 
	
	The set $R=\{s,sr,sr^2,\ldots,sr^{n-1}\} \subseteq D_{2n}$ consists of reflections. In other words, elements of order $2$. So every element $a$ of $R$ is adjacent only to $1$ the identity element and itself in the power graph. This shows that $N(a) = \{1\}$ for all $a \in R$ and the set $R$ is contained in a single twin class of the graph $G$. Since every element of $D_{2n}$ lying outside the set $R$ has the element $r$ in its neighbourhood, we have $R$ itself as a twin class. 
	

	Let $D_{2n}$ be the Dihedral group of order $2n$, where $n\geq 3$. Then $$G(D_{2n})=H[\underbrace{K_1,\dots,K_1}_{\substack{\text{corresponds to} \\ \text{reflections}}}, \underbrace{K_{\phi(n)+1}, K_{\phi(d_1)}, K_{\phi(d_2)}, \ldots, K_{\phi(d_t)}}_{\substack{\text{corresponds to} \\ \text{rotations}}}].$$
	
	
	Now, using the approach similar to the previous example, we can explicitly calculate the $m$-Steiner Wiener index of $G(D_{2n})$.
\end{example}

\begin{example} 
	Consider the group $D_{12}=\left\langle r,s:r^6=1, s^2=1, rs=sr^{-1}\right\rangle$. Then $C_1=\{1\}$, $C_2=\{r, r^5\}$, $C_3=\{r^3\}$, $C_4=\{r^2, r^4\}$ and $C_5=\{s,sr,\ldots,sr^5\}$.  Also, $G(C_i)$ is complete for $i=2,4$ and $G(C_i)=K_1$ for $i=1,3$ and $G(C_5)$ is totally disconnected.   Therefore, by Theorem \ref{pw}, $W(G(D_{12})) = (1+1)+ (0+0+ 2 (15))+ 81=113$.
\end{example}

\begin{example} 
	Consider the group $Q_8=\left\langle a, b: a^4=1, b^2=a^2, b^{-1}a b= a^{-1}\right\rangle$ Then $Q_8=\{ 1,a,a^2,a^3, b, ab, a^2b, a^3b\}$. Also, $C_1=\{1,a^2\}$, $C_2=\{a, a^3\}$, $C_3=\{ab, a^3b\}$, $C_4=\{b, a^2b\}$. Observe that $G(C_i)$ is complete for all $i\in\{1,2,3,4\}$.
	Therefore, by Theorem \ref{ppp}, $SW_6(G(Q_8)) = 27(5)+ 6= 141$.
\end{example}

\begin{example}
	Let $\Gamma$ be a finite group of order $n$ such that $o(x)=2$ for all $x\neq e$ in $\Gamma$ and let $G$ be the power graph of the group $\Gamma$. Fix $1<m<n$. 
	
	Clearly, $G \cong K_{1,n-1}$. Also $C_1=\{e\}$ and $C_2=G\smallsetminus \{e\}$. By Theorem \ref{ppp}, $$SW_m(G)= m\binom{n-1}{m}+(m-1) \binom{n-1}{m-1}.$$ 
\end{example}

\section{The Steiner Wiener index of graphs associated to finite commutative ring with unity} Let $R$ be an arbitrary finite commutative ring with unity and $I$ be an ideal of $R$.
\begin{defn}
	The zero divisor graph of $R$, denoted by $\Gamma(R)$, is defined as follows. The vertex set $V(\Gamma)$ is the set of all non-zero zero divisors of $R$. Two distinct vertices $x$ and $y$ are joined by an edge in $\Gamma$ if their product is zero.
\end{defn}
\begin{defn}
	The ideal-based zero divisor graph of $R$ ($I$-zero divisor graph in short), denoted by $\Gamma_I(R)$, is defined as follows. The vertex set $V(\Gamma_I(R))$ is the set of elements $x$ of $R\backslash I$ satisfying $xy \in I$ for some $y \in R \backslash I$. Two distinct vertices $x$ and $y$ are joined by an edge in $\Gamma_I(R)$ if their product is in $I$.
\end{defn}
\begin{defn}
	The co-maximal ideal graph of $R$, denoted by $\Gamma_{C}^I(R)$, is defined as follows. The vertex set $V(\Gamma_C^I)$ is the set of proper ideals of $R$ which are not contained in the Jacobson radical of $R$. Two vertices $I$ and $J$ are joined by an edge in $\Gamma_C^I$ if their sum is the full ring $R$.
\end{defn}






Let $G$ be any of one of the three graphs defined above. Let $V$ be the vertex set of $G$.  Define an equivalence relation $\sim$ on $V$ by $a \sim b$ if $N(a) \backslash \{b\} = N(b) \backslash \{a\}$. We observe that,

1. In $\Gamma(R)$, $a\sim b$ if $\ann(a) \backslash \{b\} = \ann(b) \backslash \{a\}$. 

2. In $\Gamma_I(R)$, $a \sim b$ if $\ann_I(a) \backslash \{b\} = \ann_I(b) \backslash \{a\}$ 

where $\ann_I(a) = \{r \in R-I : ra \in I\}$.

3. In $\Gamma_C^I(R),  I \sim J$ if $\{M : I + M = R\}  \backslash \{J\} = \{M: J + M = R\} \backslash \{I\}$.


Let $G$ be one of the above three graphs with the respective equivalence relation defined above.  Let $C_1, C_2, \dots, C_k$ be the equivalence classes (twin classes) of this equivalence relation with respective representatives $c_1,\dots,c_k$. Then by Proposition \ref{join}, we have
$$G = H[G[C_1],\dots,G[C_k]]$$ where $H = G[\{c_1,\dots,c_k\}]$ is the subgraph of $G$ induced by the set of representatives. Hence, we can apply Theorem \ref{mainthm}. We get

\begin{thm}\label{id5}
	 Let $R$ be a finite commutative ring with unity. Let $I$ be an ideal in $R$. Then  the  $m$-Steiner Wiener index of $\Gamma_I(R)$ is given by

 	{\small \begin{equation}
		SW_m(\Gamma_I(R)) = \sum_{\substack{1 \le i \le k \\ \Gamma_I(C_i) \text{ is complete}}}(m-1) \binom{n_i}{m} + \sum_{\substack{1 \le i \le k \\ \Gamma_I(C_i) \text{ is empty}}}m \binom{n_i}{m}+\sum_{\substack{S \in \mathcal P_m(V) \\ |I(S)|=l>1}} \Big(d_H^l(S_r)+ \sum\limits_{c_i \in S_r}(t_i - 1)\Big).
	\end{equation}} where $|C_i| = n_i$
\end{thm}

We explain our result with the following working examples.
\begin{example}
Consider the ring $R=\mathbb{Z}_{24}$ and $I=(8)$. Then $V(\Gamma_I(R))=\{2,4,6,10,12,14,\\ 18,20,22\}$. Also, $C_1=\{2,6, 10, 14, 18, 22\}$ and $C_2=\{4,12,20\}$  are the twin classes of $\Gamma_I(R)$. The graphs  $\Gamma_I(C_2)$ is complete and $\Gamma_I(C_1)$ is totally disconnected. Clearly, $H\cong K_{2}$. For each $i\in\{2,4,6,10,12,14,18,20,22\}$, define 
$S^i=V(\Gamma_I(R))-\{i\}$ and $|S^i|=8$. Then $S_r^i=\{c_1,c_2\}$ and $d_H^l(S_r^i)=2$ for all $i\in\{2,4,6,10,12,14,18,20,22\}$.  By using Theorem  \ref{id5}, $SW_8(\Gamma_C^I(R))=9(2+5)=63$. 

\end{example}

 \begin{example}
Consider the ring $R=\frac{\mathbb{Z}_2[x]}{(x^3)}\times\mathbb{Z}_2$ and $I=(0)\times\mathbb{Z}_2$. Then $V(\Gamma_I(R))=\{u_1=(\overline{x},1), u_2=(\overline{x^2},0), u_3=(\overline{x+x^2},1), u_4=(\overline{x},0), u_5=(\overline{x^2},1), u_6=(\overline{x+x^2},0)\}$. 

Suppose $m=4$. Consider $S^1=\{u_1,u_2,u_3,u_4\}$, $S^2=\{u_1,u_2,u_3,u_5\}$, $S^3=\{u_1,u_2,u_3,u_6\}$, $S^4=\{u_1,u_2,u_4,u_5\}$, 
$S^5=\{u_1,u_2,u_4,u_6\}$, $S^6=\{u_1,u_3,u_4,u_5\}$, $S^7=\{u_1,u_3,u_4,u_6\}$, $S^8=\{u_2,u_3,u_4,u_5\}$,  $S^9=\{u_3,u_2,u_4,u_6\}$, $S^{10}=\{u_3,u_5,u_4,u_6\}$, $S^{11}=\{u_5,u_3,u_2,u_6\}$, $S^{12}=\{u_2,u_5,u_4,u_6\}$, $S^{13}=\{u_1,u_5,u_4,u_6\}$, $S^{14}=\{u_1,u_2,u_5,u_6\}$, $S^{15}=\{u_1,u_3,u_5,u_6\}$. Then $d(S^7)=4$,$d(S^i)=3$ for all $i\neq 7$ and hence by definition \ref{sw2}, $SW_4(\Gamma_I(R))=46$.

Also, $C_1=\{u_2,u_5\}$, $C_2=\{u_1,u_3,u_4,u_6\}$ are twin classes of $\Gamma_I(R)$ and $H\cong K_2$.  Clearly $S_r^7=\{c_2\}$ and $S_r^i=\{c_1,c_2\}$ for all $i\in\{1,2,\ldots,6, 8,9,\ldots, 15\}$. Therefore by Theorem \ref{id5}, $SW_4(\Gamma_I(R))=46$.

\end{example} 


\begin{lem}\label{id6}
Let $R$ be a finite commutative ring with unity. The the  Wiener index of $\Gamma_I(R)$ is given by 
\begin{equation}
		W(\Gamma_I(R)) = \sum_{\substack{1 \le i \le k \\ \Gamma_I(C_i) \text{ is complete}}}\binom{n_i}{2} + \sum_{\substack{1 \le i \le k \\ \Gamma_I(C_i) \text{ is empty}}}2\binom{n_i}{2}+\sum_{1 \le i<j \le k}n_in_jd_H(i,j).
	\end{equation} where $|C_i| = n_i$.  
\end{lem}

 \begin{example}
Consider the ring $R=\mathbb{Z}_6\times\mathbb{Z}_2$ and $I=(0)\times \mathbb{Z}_2$. Then $V(\Gamma_I(R))=\{(2,0),(3,0),(2,1),(3,1), (4,0),(4,1)\}$. Also, $C_1=\{(3,0),(3,1)\}$,  $C_2=\{(2,0), (2,1),  (4,0),\\(4,1)\}$, $\Gamma_I(C_i)$ is totally disconnected and so $H(R)\cong K_2$. By using Lemma \ref{id6}, $W(\Gamma_I(R))=2(1)+ 2(6)+2(4)(1)=22$. 

 \end{example}


In the following, the comaximal ideal graph $\Gamma_C^I(R)$ of an Artinian commutative nonlocal ring $R$ with unity is realized as a generalized composition of suitable choices of graphs. We prove some graph theoretic properties of these classes. 

\begin{lem}\label{a h connected}
Let $R$ be a commutative Artinian nonlocal ring with unity. Then $H(R)$ is connected.
\end{lem}
\begin{pf}
	Let $c_i$ and $c_k$ be two distinct vertices of $H$. Then $N(c_i) \ne N(c_k)$. Let $x_1 \in N(c_i)$ and $x_2 \in N(c_k)$. Suppose $x_1+x_2\neq R$. Then $c_i \sim x_1+x_2 \sim c_k$ is a path joining $c_i$ and $c_k$. Suppose $x_1+x_2=R$. Then  then $c_i \sim x_1 \sim x_2 \sim c_k$ is a path joining $c_i$ and $c_k$
\end{pf}

In the following proposition, we prove that the graph $\Gamma_C^I(R)$ is a generalized composition
of certain empty graphs.  This result is due to Theorem \ref{mainthm}.

\begin{prop}\label{c join}
	Let $G_i$ be the subgraph induced by the set $C_i$ in $\Gamma_C^I(R)$. Then $\Gamma_C^I(R) = H[G_1,\dots,G_t]$ and $G_i$ is totally disconnected for $1\leq i\leq k$. 
\end{prop}

 Let $R$ be an  Artinian commutative nonlocal ring with unity.  In this
subsection,  by using Theorem \ref{mainthm}, we deduce a formula for the $m$-Steiner Wiener index of $\Gamma_C^I(R)$. Then by Proposition \ref{c join},
$\Gamma_C^I(R) = H[G_1,\dots,G_t]$ and $G_i$ is totally disconnected for $1\leq i\leq t$ and also by Lemma \ref{general distance}, $S_r^c=\emptyset$ and $S_r^e=S_r$. Hence, we can use Theorem \ref{mainthm} to calculate the $m$-Steiner Wiener index of $\Gamma_C^I(R)$.

\begin{thm}
	 Let $R$ be an  Artinian commutative nonlocal ring with unity. With the notations as in Section \ref{ks}, the  $m$-Steiner Wiener index of $\Gamma_C^I(R)$ is given by

 	\begin{equation}
		SW_m(\Gamma_C^I(R)) =  \sum_{\substack{1 \le i \le k \\ G_i \text{ is empty}}}m \binom{n_i}{m}+\sum_{\substack{S \in \mathcal P_m(V) \\ |I(S)|=l>1}} \Big(d_H^l(S_r)+  +\sum\limits_{c_i \in S_r}(t_i-1)\Big).
	\end{equation} where $|G_i|= n_i$. 
\end{thm}

\begin{example}
Consider the ring $R=\mathbb{Z}_2\times\mathbb{Z}_2\times\mathbb{Z}_4$. Then $V(\Gamma_C^I(R))=\{u_1,u_2,\ldots,u_9\}$, where $u_1=\mathbb{Z}_2\times\mathbb{Z}_2\times(2)$, $u_2=(0)\times(0)\times \mathbb{Z}_4$, $u_3=\mathbb{Z}_2\times\mathbb{Z}_2\times(0)$, $u_4=\mathbb{Z}_2\times(0)\times\mathbb{Z}_4$, $u_5=(0)\times\mathbb{Z}_2\times(0)$, $u_6=(0)\times \mathbb{Z}_2\times (2)$, $u_7=(0)\times\mathbb{Z}_2\times\mathbb{Z}_4$,
$u_8=\mathbb{Z}_2\times(0)\times(0)$, $u_9=\mathbb{Z}_2\times(0)\times(2)$. 

For each $i=1$ to $9$, define 
$S^i=V(\Gamma_C^I(R))-\{u_i\}$ and $|S^i|=8$. Then $d(S^4)=d(S^7)=8$ and $d(S^i)=7$ for all $i\notin\{4,7\}$. Hence by definition \ref{sw2}, $SW_8(\Gamma_C^I(R)=65$.

Also $C_1=\{u_1,u_3\}$,  $C_2=\{u_7\}$, $C_3=\{u_4\}$, $C_4=\{u_2\}$, 
$C_5=\{u_5,u_6\}$, $C_6=\{u_8,u_9\}$ are twin classes of $\Gamma_C^I(R)$. $S_r^2=\{c_1,c_2,c_3,c_5,c_6\}$, $S_r^4=\{c_1,c_2,c_4,c_5,c_6\}$, $S_r^7=\{c_1,c_3,c_4,c_5,c_6\}$, and $S_r^i=\{c_1,\ldots,c_6\}$ for $i\in\{1,3,5,6, 8,9\}$. Also $d_H^l(S_r^2)=4$ and $d_H^l(S_r^i)=5$ for $i\in\{1,3,4,\ldots, 9\}$.  By using Theorem \ref{mainthm}, $SW_8(\Gamma_C^I(R))=65$.
\end{example} 

By using Theorem \ref{mainthm}, we deduce a formula for the Wiener index of $\Gamma_C^I(R)$.

\begin{thm}\label{www}
	 Let $R$ be an  Artinian commutative nonlocal ring with unity. Then the  Wiener index of $\Gamma_C^I(R)$ is given by

 	\begin{equation}
				W(\Gamma_C^I(R)) =  \sum_{\substack{1 \le i \le t \\ G_i \text{ is empty}}}2\binom{n_i}{2}+\sum_{1 \le i<j \le t}n_in_jd_H(i,j).
	\end{equation} where $|C_i|= n_i$. 
\end{thm}

Next, we explain the above theorem with various examples.

\begin{example} 
	Consider the ring $R = \mathbb Z_{p_1^3p_2^2} \cong \mathbb Z_{p_1^3} \times \mathbb Z_{p_2^2}$ where $p_1$ and $p_2$ are distinct prime numbers. We will calculate the  Wiener index $W(\Gamma_C^I(R))$ using Theorem \ref{www}. We have $n=2$, $C_1=\{(p_1)\times\mathbb Z_{p_2^2}, (p_1^2)\times\mathbb Z_{p_2^2}, (0)\times\mathbb Z_{p_2^2}\}$, $C_2=\{\mathbb Z_{p_1^3}\times (0), \mathbb Z_{p_1^3}\times (p_2)\}$ are twin classes of $\Gamma_C^I(R)$ and $H\cong K_2$.  Therefore, by Theorem \ref{www}, $W(\Gamma_C^I(R)) = 2 (1)+ 2 (3)+2 \cdot 3 (1)=14.$
\end{example}

\begin{example} 
Consider the ring $R=\mathbb{Z}_3\times\mathbb{Z}_5\times\mathbb{Z}_9$. Then $V(\Gamma_C^I(R))=\{v_1,v_2,\ldots,v_9\}$, where $v_1=\mathbb{Z}_3\times\mathbb{Z}_5\times(3)$, $v_2=(0)\times(0)\times \mathbb{Z}_9$, $v_3=\mathbb{Z}_3\times\mathbb{Z}_5\times(0)$, $v_4=\mathbb{Z}_3\times(0)\times\mathbb{Z}_9$, $v_5=(0)\times\mathbb{Z}_5\times(0)$, $v_6=(0)\times \mathbb{Z}_5\times (3)$, $v_7=(0)\times\mathbb{Z}_5\times\mathbb{Z}_9$,
$v_8=\mathbb{Z}_3\times(0)\times(0)$, $v_9=\mathbb{Z}_3\times(0)\times(3)$. Note that $C_1=\{v_1,v_3\}$,  $C_2=\{v_7\}$, $C_3=\{v_4\}$, $C_4=\{v_2\}$, 
$C_5=\{v_5,v_6\}$, $C_6=\{v_8,v_9\}$ are twin classes of $\Gamma_C^I(R)$. Therefore, by Theorem \ref{www},
$W(\Gamma_C^I(R)) = 2 (1)+ 2 (0)+2 (0)+2 (0)+2(1) +2(1) + 22+9+8+12+12= 69$.
\end{example}

From the above examples, we see that the Wiener index and the more general Steiner Wiener index of graphs can be explicitly calculated using our formula. This shows the applicability of our result.



\begin{thebibliography}{99}
	
\bibitem{abs22}  Ali Azimi, Ravindra B. Bapat and Shivani Goel, Steiner distance matrix of caterpillar graphs, Special Matrices,  10 (1)  267--284, 2022.

\bibitem{aja2019} Abdussakir, Sudarman, M N Jauhari and F Ali, Survey on topological indices and graphs associated with a commutative ring, J. Phys.: Conf. Ser. 1562 (2020) 012008.
	
\bibitem{al99} D. F. Anderson, P. S. Livingston, The zero-divisor graph of a commutative ring, J. Algebra  217 (1999)  434--447.

\bibitem{survey1} Andrey A. Dobrynin, Roger Entringer, Ivan Gutman,  Wiener Index of Trees: Theory and Applications, Acta Applicandae Math.   66 (2001) 211--249.





\bibitem{bst93} F. Boesch, C. Suffel, R. Tindell, The neighborhood inclusion structure of a graph, Math. Comput. Modelling 17 (1993) 25–28.

\bibitem{cmica10} Carmen Hernando, Mercè Mora, Ignacio M. Pelayo, Carlos Seara, David R. Wood, Extremal Graph Theory for Metric Dimension and Diameter, the Electronic Journal of Combinatorics,  2010(17), \#R30.



\bibitem{d20}  Daniel Weißauer, Isometric subgraphs for Steiner distance, Journal of Graph Theory,  94 (4)  597--613, 2020.

\bibitem{da08} David Auger, Induced Paths in Twin-Free Graphs, the Electronic Journal of Combinatorics,  2008(15), \#N17.

\bibitem{dsb} Deepa Sinha, Bableen Kaur,  Some aspects of zero-divisor graphs for the ring of Gaussian integers modulo $2^n$, J. Appl. Math. Comput. 68 (2022) 69-81

\bibitem{dbw} Douglas. B. West, Introduction to Graph Theory, 2nd Ed., Pearson Education, Singapore  2001.

\bibitem{dp20} Damaschke, Peter, Equitable induced decompositions of twin graphs, Australas. J. Combin. 76 (2020), part 1, 24–40.


\bibitem{fh00} D.Frank Hsu, X.-D.Hu, On shortest three-edge-connected Steiner networks with Euclidean distance, Discrete. Appl. Math.   103, Issues 1–3, 141-152, 2000.

\bibitem{gb21} Giuliano Basso, False-twin-free graphs with a fixed number of negative eigenvalues, Linear algebra and its Applications, 618(2021), 144-149.


\bibitem{hh86} Hans-Jürgen Bandelt, Henry Martyn Mulder, Distance-hereditary graphs, Journal of Combinatorial Theory, Series B, 41(2), 1986, 182-208.



\bibitem{ioa15} Iiro Honkalaa, Olivier Hudry, Antoine Lobsteinc, On the number of optimal identifying codes in a twin-free graph, Discrete. Appl. Math.   2015 (180), 111-119.

\bibitem{iioa07} Irène Charon, Iiro Honkala, Olivier Hudry, Antoine Lobstein, Structural Properties of Twin-Free Graphs, the Electronic Journal of Combinatorics,  2007(14), \#R16.




\bibitem{lo08} Linda Eroh, Ortrud R.Oellermann, Geodetic and Steiner geodetic sets in 3-Steiner distance hereditary graphs, Discrete. Math.   308(18), 4212-4220, 2008.




\bibitem{mrt14} Maria Axenovich, Ryan Martin, Torsten Ueckerdt, Twins in graphs, European Journal of Combinatorics, 2014(39), 188-197.

\bibitem{mk09}  Matthias Kriesell, Edge disjoint Steiner trees in graphs without large bridges, Journal of Graph Theory,  62 (2)  188--198, 2009.

\bibitem{oo99} On Steiner centers and Steiner medians of graphs, Networks, 34(4), 258-263, 1999.

\bibitem{pho97} Peter Dankelmann, Henda C.Swart, Ortrud R.Oellermann, On the average Steiner distance of graphs with prescribed properties, Discrete. Appl. Math.   79, Issues 1–3, 91-103, 1997.






\bibitem{sch} A.J. Schwenk,  Computing the characteristic polynomial of a graph, in: R. Bary, F. Harary (Eds.), Graphs Combinatorics, in: Lecture Notes in Mathematics, Springer-Verlag, Berlin,	406 (1974) 153--172.

\bibitem{sga} K. Selvakumar, P. Gangaeswari, G. Arunkumar,  The Wiener index of the zero-divisor graph of a finite commutative ring with unity, Discrete Applied Mathematics 311 (2022) 72–84.

\bibitem{tfz08} B.-S. Tam, Y.-Z. Fan, J. Zhou, Unoriented Laplacian maximizing graphs are degree maximal, Linear Algebra Appl. 429 (2008)
735–758.


\bibitem{w} H. Wiener,  Structural determination of paraffin boiling points, J. Am. Chem. Soc.  69 (1) (1947)  17--20. 

\bibitem{woh94} Wayne Goddard, Ortrud R.Oellermann, Henda C.Swart, Steiner distance stable graphs, Discrete. Math. 132(1-3), 65-73, 1994.

   





	
	

	
\end{thebibliography}
\end{document}